\documentclass{article}
\usepackage{amsmath}
\usepackage[psamsfonts]{amssymb}
\usepackage{cmmib57}
\newcommand{\bPf}{\par\vspace*{-4pt}\indent{\sc Proof.}\enskip}

\def\QED{\hskip0.1em\hfill\null\ \null\nobreak\hfill\kern3pt\vbox{\hrule\hbox
   {\vrule\kern1pt\vbox{\kern1.7pt\hbox{$\scriptscriptstyle{QED}$}
    \kern0.2pt}\kern1pt\vrule}\hrule}}

\def\END{\hskip0.1em\hfill\null\ \null\nobreak\hfill\kern3pt\vbox{\hrule\hbox
   {\vrule\kern1pt\vbox{\kern1.7pt\hbox{$\,\,\,\vspace{5pt}$}
    \kern0.2pt}\kern1pt\vrule}\hrule}}
\newtheorem{theorem}{Theorem}
\newtheorem{lemma}{Lemma}
\newtheorem{corollary}{Corollary}
\newtheorem{proposition}{Proposition}
\newtheorem{remark}{Remark}
\newtheorem{definition}{Definition}
\newtheorem{example}{Example}
\newcommand{\bCd}{\bEq\begin{CD}}
\newcommand{\eCd}{\end{CD}\eEq}
\newcommand{\bcd}{\beq\begin{CD}}
\newcommand{\ecd}{\end{CD}\eeq}
\newcommand{\ben}{\begin{enumerate}}
\newcommand{\een}{\end{enumerate}}
\newcommand{\bEq}{\begin{eqnarray}}
\newcommand{\eEq}{\end{eqnarray}}
\newcommand{\beq}{\begin{eqnarray*}}
\newcommand{\eeq}{\end{eqnarray*}}
\newcommand{\bDf}{\begin{definition}\em}
\newcommand{\eDf}{\end{definition}}
\newcommand{\bLm}{\begin{lemma}}
\newcommand{\eLm}{\end{lemma}}
\newcommand{\bPr}{\begin{proposition}}
\newcommand{\ePr}{\end{proposition}}
\newcommand{\bTh}{\begin{theorem}}
\newcommand{\eTh}{\end{theorem}}
\newcommand{\bCr}{\begin{corollary}}
\newcommand{\eCr}{\end{corollary}}
\newcommand{\bRm}{\begin{remark}\em}
\newcommand{\eRm}{\end{remark}}
\newcommand{\bEx}{\begin{example}\em}
\newcommand{\eEx}{\end{example}}

\newcommand{\ie}{{\em i.e$.$} }

\newcommand{\A}{\forall}

\newcommand{\mto}{\mapsto}

\newcommand{\cD}{\mathcal{D}}

\newcommand{\sub}{\subset}

\newcommand{\alp}{\alpha}
\newcommand{\bet}{\beta}

\newcommand{\tht}{\theta}

\newcommand{\lam}{\lambda}
\newcommand{\sig}{\sigma}

\newcommand{\ome}{\omega}

\newcommand{\Ome}{\Omega}

\title{\large{\bf On the Geometry of B\"acklund Transformations}}
\author{{\normalsize M.
Palese and E. Winterroth}\thanks{Both of them supported by
GNFM of INdAM and University of Torino proj. {\em Giovani Ricercatori 2001}. E.W.
partially supported also by University of
Erlangen--N\"urnberg and CNR grant n. 201.21-00.01.01.}
\\{\footnotesize Department of Mathematics,
University of Torino}
\\{\footnotesize Via C. Alberto 10, 10123 Torino, Italy}\\ 
{\footnotesize e--mails: 
{\sc palese@dm.unito.it, ekkehart@dm.unito.it}}}
\date{}
\begin{document}

\maketitle

\begin{abstract}

The geometry of an admissible 
B\"acklund transformation for an exterior differential system is described by
an admissible Cartan connection for a geometric structure on 
a tower with infinite--dimensional skeleton which is the universal prolongation of a
$|1|$--graded semi-simple Lie algebra.

\medskip


\end{abstract}

\section{Towers on skeletons}\label{2}

In the following we consider
infinite--dimensional objects in
the $\mathfrak{p}$--category of objects obtained as {\em projective limits} of
finite--dimensional ones \cite{Mo93}.

\bDf An
algebraic skeleton on a finite--dimensional vector space $V$ is a triple
$(E,G,\rho)$, with $G$ a
$\mathfrak{p}$--Lie group, $E=V\oplus\mathfrak{g}$, $\mathfrak{g}$ the Lie
algebra of $G$, and $\rho$ a representation of $G$ on $E$ such that
$\rho(g)x=Ad(g)x$, for $g\in G$, $x\in\mathfrak{g}$. 
An infinitesimal skeleton can be analogously defined {\em via} the
representation of $\mathfrak{g}$ on $E$.\END \eDf

\bDf Let $(E,G,\rho)$ be a skeleton on $V$ and $Z$ a manifold of type
$V$ \cite{Mo93}. We say that a $\mathfrak{p}$--principal fibre bundle
$P(Z,G)$ provided with an absolute parallelism $\omega$ on $P$ is a {\em
tower} on $Z$ with skeleton $(E,G,\rho)$ if $\omega$ takes values in $E$
and satisfies: $R^{*}_{g}\omega = \rho(g)^{-1}\ome$, for
$g\in G$; $\ome(\tilde{A}) = A$, for $A\in\mathfrak{g}$;
here $R_{g}$ denotes the right translation and $\tilde{A}$ the fundamental
vector field induced on $P$ from $A$.\END\eDf

\subsection{Cartan connections}

Let $\mathfrak{g}$ be a Lie algebra and $\mathfrak{k}$ a Lie subalgebra of
$\mathfrak{g}$. Let $K$ be a Lie group with Lie algebra $\mathfrak{k}$ equipped
with a representation $Ad: K\to GL(\mathfrak{g})$ such that its differential
coincides with the adjoint representation of $\mathfrak{k}$ on $\mathfrak{g}$.

\bDf Let $P(Z,K)$ be a principal fibre bundle over a manifold $Z$ with
structure group $K$. A {\em Cartan connection} in $P$ of type
$(\mathfrak{g},K)$ is a $1$--form $\omega$ on $P$ with values in
$\mathfrak{g}$ satisfying the following conditions: 
$\omega|_{T_{u}P}: T_{u}P\to \mathfrak{g}$ is an isomorphism $\forall u\in
P$; $R^{*}_{g}\omega=Ad(g)^{-1}\omega$ for $g\in K$;
$\omega(\tilde{A})=A$ for $A\in \mathfrak{k}$.\END\eDf
$(\mathfrak{g},{K},Ad)$ is a skeleton on $V$, with
$\mathfrak{g}=\mathfrak{k}\oplus V$. Then it is clear that a Cartan connection
$(P,Z,K,\omega)$ of type $(\mathfrak{g},K)$ is a tower on $Z$.

In the following we assume the Lie algebra $\mathfrak{g}$ to be a generalized
semi-simple $|1|$-graded Lie algebra \ie
$\mathfrak{g}=\mathfrak{g}_{-1}\oplus\mathfrak{g}_{0}\oplus
\mathfrak{g}_{1}$ \cite{Oc70}.

\bRm
According with the Lie algebra $|1|$--grading the Cartan connection 
form $\ome$ and 
its curvature $\kappa$ split as $\ome= \ome_{-1}\oplus\ome_{0}\oplus\ome_{1}$ and 
$\kappa= \kappa_{-1}\oplus\kappa_{0}\oplus\kappa_{1}$.\END\eRm

\bDf Let $G$ be a semi-simple Lie group, with $|1|$--graded Lie algebra
$\mathfrak{g}$ as above and $K$ the closed subgroup of $G$ corresponding to
the Lie algebra $\mathfrak{g}_{0}\oplus \mathfrak{g}_{1}$. 
A {\em $K$--structure} on
$Z$ is a principal fiber bundle $P\to Z$ with structure group $K$ equipped
with a {\em soldering} form $\tht=\tht_{-1}\oplus\tht_{0}\in
\Ome^{1}(P,\mathfrak{g}_{-1}\oplus\mathfrak{g}_{0})$ such that: $\tht_{-1}(\xi)=0$, if and only if $\xi$ is a vertical
vector; $\tht_{0}(\tilde{X}+\tilde{Z})=Y$, $\A Y\in\mathfrak{g}_{0},
Z\in\mathfrak{g}_{1}$; $(R_{b})^{*}\tht =Ad(b^{-1})\tht$, $\A b\in K$,
where $Ad$ means the action on the vector space
$\mathfrak{g}_{-1}\oplus\mathfrak{g}_{0}\simeq\mathfrak{g}/\mathfrak{g}_{1}$
induced by the adjoint action.\END\eDf

\bDf Let $(P,\tht)$ be a $K$--structure on $Z$. A Cartan connection $\ome$
on $P$ is called {\em admissible} if and only if it is of the form
$\ome=\tht_{-1}	\oplus\tht_{0}\oplus\ome_{1}$.\END\eDf

\bRm\label{torsion2} If the $K$--structure has zero torsion, \ie if it is a
reduction of $L^{2}(Z)\to Z$ to $K$, then the curvature of the induced
Cartan connection $\ome=\tht_{-1}\oplus\tht_{0}\oplus\ome_{1}$ is such that
$\kappa_{-1}=0$ \cite{CSS94}.\END\eRm

\section{B\"acklund transformations and induced Cartan connections}\label{3}

Let $\pi:  U\to  X$, $\tau: Z\to  X$,
be two (vector) bundles 
with local fibered coordinates $(x^{\alp},u^{A})$ and $(x^{\alp},z^{i})$,
respectively, where $\alp=1,\ldots,m=\textstyle{dim} X$, $A=1,\ldots,n=\textstyle{dim} U -
\textstyle{dim} X$, $i=1,\ldots,N=\textstyle{dim}Z-\textstyle{dim} X$.  
A system of nonlinear field equations of order $k$ on $ U$ is geometrically described as an
exterior differential system $\nu$ on $J^{k} U$. The solutions of the field
equations are (local)  sections $\sig$
of $U\to X$ such that $(j^{k}\sig)^{*}\nu=0$. We shall also denote by
$J^{\infty}\nu$ ({\em resp.} $j^{\infty}\sig$) 
the infinite order jet prolongation of $\nu$ ({\em resp.} $\sig$). 

\subsection{Admissible B\"acklund transformations}

Let $B$ be the infinite--order contact transformations group on 
$J^{\infty}U$.

\bDf\label{admissible}
The group of (infinitesimal) {\em B\"acklund transformations for the
system $\nu$} is the closed subgroup
${ \tilde K}$ of $ B$ which leaves invariant solution submanifolds 
of $J^{\infty}\nu$.
The group of (infinitesimal) {\em generalized B\"acklund transformations for the system
$\nu$} is the closed subgroup $K$ of
$B$ which leaves invariant $J^{\infty}\nu$ \cite{AI79}.\END\eDf

Let $\pi: U\to X$, $\tau:Z\to X$, be vector bundles as the above and 
$\pi^{1}: {J^{1} U} \to X$, $\tau^{1}: {J^{1}Z} \to X$, the first order jet
prolongations bundles, with local fibered coordinates
$(x^{\alp},u^A,u^{A}_{\alp})$,
$(x^{\alp},z^i,z^{i}_{\alp})$, respectively. Furthermore, let
$(\partial_{\beta}$, $\partial_{A}$, $\partial_{A}^{\beta})$,
$(\partial_{\beta},\partial_{i}$, $\partial_{i}^{\beta})$ and $(dx^{\beta}$,
$du^{A}$, $du^{A}_{\beta})$, $(dx^{\beta},dz^{i}$, $dz^{i}_{\beta})$ be local
bases of tangent vector fields and $1$--forms on $J^{1} U$ and $J^{1}Z$,
respectively.

\bDf
We define a B\"acklund map to be 
the fibered morphism over $Z$:
$\phi: J^{1} U\times_{X}Z\to J^{1}Z: (x^{\alp},u^A,u^{A}_{\alp};z^{i})\mto
(x^{\alp},z^{i},z^{i}_{\alp})$,
with $z^{i}_{\alp}=\phi^{i}_{\alp}(x^{\bet},u^A,u^{A}_{\bet};z^{j})$.

The fibered morphism $\phi$ is said to be an {\em admissible} B\"acklund
transformation for the differential system $\nu$ if
$\phi^{i}_{\alpha}=\cD_{\alpha}\phi^{i}$ and the integrability
conditions coincide with the
exterior differential system $\nu$. 
\END\eDf

\bRm 
By pull--back of the contact structure on $J^1 Z$, the B\"acklund morphism 
induces an horizontal distribution, the
{\em induced B\"acklund connection},
on the bundle $(J^{1}{ U}\times_{X}{Z},$ $J^{1}{ U},$
$\pi^{1*}_{0}(\eta))$ \cite{PRS79}.\END\eRm

\bTh\label{maintheorem}
The following statements are equivalent \cite{PaWi02} .
\begin{enumerate}
\item $\phi$ is an admissible B\"acklund transformation for the
differential system $\nu$. 
\item The induced B\"acklund connection is
$\bar{K}$--invariant, where $\bar{K}$ is a normal subgroup
$\bar{K}\sub ({\tilde K} \cap K) \sub  B$ leaving invariant (the infinite order
prolongation of) $\nu$ and its solutions.
\end{enumerate}
\eTh

Let now $Z$ be a vector bundle over the basis $X$, the fibers of which are modelled over
the homogeneous space $\tilde{K}/\bar{K}$ such that fibers are vector spaces
of type $\bar{\mathfrak{K}}_{-}=\oplus_{p<0}\,\,\bar{\mathfrak{K}}_{p}$, with
$\bar{\mathfrak{K}}_{-}$ a graded abelian Lie algebra. Let
$P$ be a tower on $Z$ with algebraic skeleton
$(\tilde{\mathfrak{k}},\bar{\mathfrak{K}},Ad)$, where $\tilde{\mathfrak{k}}$ is the Lie
algebra of $\tilde{K}$. Suppose
$U$ be a vector bundle (over the same basis $X$) with a left action $\lam$ of $\bar{K}$ on
$U$ (as a manifold). For each tower
$(P,Z,\bar{K},\omega)$ we have a vector bundle
$U_{\lam}(Z)=P\times_{\bar{K}} U$ over $X$ and {\em vice versa}.
 
Assume $\tilde{\mathfrak{k}}$ to be the universal prolongation \cite{Ta79} of a
$|1|$--graded semi-simple Lie algebra $\mathfrak{g}$ such that $\bar{\mathfrak{k}}$ and
$\bar{\mathfrak{K}}_{-}$ are (the prolongation of) 
$\mathfrak{g}_{0}\oplus\mathfrak{g}_{1}$ and $\mathfrak{g}_{-1}$, respectively.

\bTh\label{main}
A B\"acklund transformation admissible for an exterior
differential system induces the tower $(P,Z,\bar{K},\omega)$, where $\ome$ is an
admissible Cartan connection for a $\bar{K}$-structure over $Z$. 

\eTh

\bPf 
It follows from Remark 
\ref{torsion2} and Theorem \ref{maintheorem}. 
In fact, a B\"acklund morphism can be seen as a 
$\bar{K}$--equivariant section of the
bundle $ U_{\lam}(Z) \to Z$; it induces a reduction of
$L^{2}(Z)\to Z$ defining a
$\bar{K}$--structure on $Z$ with zero torsion. 
The admissible
Cartan connection is the tower $(P,\bar{K},\tht)$ induced \cite{CSS94} from a 
$\mathfrak{g}_{0}$--principal connection on the underlying first order structure.
\QED

As a consequence one can built a cohomological theory of complete integrability 
for (nonlinear) exterior differential systems.
Cohomological conditions are in fact given for a graded simple Lie algebra to be a
universal prolongation \cite{Ta79}. This topic will be developed elsewhere.

\section*{Acknowledgments}
Thanks are due to M. Francaviglia, R.A. Leo and G. Soliani for useful discussions.

 
\end{document}